\documentclass{article}
\newcommand\G{\mathbf{G}}
\newcommand\T{\mathbf{T}}
\newcommand\D{\mathbf{D}}
\newcommand\N{\mathbf{N}}

\begin{document}
\title{On the number of open knight's tours}
\author{
H\'ector Cancela\footnote{
Facultad de Ingenier\'ia,
Instituto de Computaci\'on.
J. Herrera y Reissig 565. CP 11300,
Montevideo. Uruguay.
cacela@fing.edu.uy}
\qquad
Ernesto Mordecki\footnote{
Facultad de Ciencias,
Centro de Matem\'atica.
Igu\'a 4225. CP 11400,
Montevideo. Uruguay.
mordecki@cmat.edu.uy}
}
\date{\today}
\maketitle
\begin{abstract}
We review the state of the art in the problem of counting the number \emph{open knight tours}, since the publication in internet of a computation of this quantity \cite{ab}. 
\end{abstract}


\section{Problem description}

``A knight's tour is a sequence of moves of a knight on a chessboard such that the knight visits every square only once. If the knight ends on a square that is one knight's move from the beginning square (so that it could tour the board again immediately, following the same path), the tour is closed, otherwise it is open.
The knight's tour problem is the mathematical problem of finding a knight's tour"\cite{wp}.
In this note we are interested in the determination of the \emph{number of solutions}
of the knight's tour problem, i.e. in the computation of the number of open knight tours. 

The history of the problem dates back to the Hindu culture, 
in the 9th century BC \cite{wp}. 
Many mathematicians made contributions to this problem
finding solutions with special properties, and providing methods to find them.
The list includes De Moivre, Euler, Legendre and
Vandermonde \cite{guik,biggs,p,w}.

To define precisely the quantities we are interested in, we borrow the notation from 
George Jelliss's web page \cite{gj}, denoting by
\begin{description}
\item[$\G$] the number of $\G$eometrically distinct open tours.
\item
[$\T$] the number of open $\T$our diagrams, by rotation and symmetry $\T = 8\,\G$.
\item[$\N$] the number of open tour $\N$umberings: 
$\N = 2\,\T = 16\,\G$ since each can tour diagram be numbered from either end.
\item[$\D$]
 is the number of closed tour $\D$iagrams.
\end{description} 
These quantities refer to the standard $8\times 8$ chessboard.

\section{Number of closed tours}

The number of closed tours was independently computed by McKay \cite{mk} and 
Wegener \cite{wegener}, obtaining that the number of classes of re-entrant solutions
(not taking into account the order and the initial square) is
$\D=13,267,364,410,532$, that is approximately $10^{13}$.

\section{An estimation of the number of open tours}

The number of open tours in an a classical $8\times 8$ chessboard remained 
an open question, as reported in a comprehensive web page about Knight tours \cite{gj}, where some bounds and estimation are discussed. The method of importance sampling combined with the Warnsdorff's Rule \cite{guik}
was used by Cancela and Mordecki in 2006 \cite{cm} to conjecture that 
$G\simeq1.22 \times 10^{15}$, providing a $99\%$-confidence interval of the form
$[1.220,1.225]\times 10^{15}$.

\section{A computation of the number of open tours}

A comprehensive computation of the total number of open tours was carried out by Alexander Chernov, and published on the Internet \cite{ab}. The final amount is   
$\T=9,795,914,085,489,952$ open tours diagrams, that gives
$\G=1,224,489,260,686,244$.
This number belongs to the confidence interval in \cite{cm}.
This number is also published as the 8-th term of the sequence A165134 in \cite{seq}, by A. Chernov. Although we believe that Chernov's result is correct, we expect to have an independent verification, that we hope will be carried out soon.


\begin{thebibliography}{99}

\bibitem{ab}
Alex Chernov.
O kolichestve nezamknutykh marshrutov konya (Russian). 
On the number of open knight tours  (05/10/2014)
Retrieved July 6, 2015.
http://alex-black.ru/article.php?content=141.

\bibitem{cm}
H. Cancela and E. Mordecki.
Counting Knight's Tours through the Randomized Warnsdorff Rule
arXiv:math/0609009. (2006)
\bibitem{guik}
E. Ya. Guik.
{\it Chess and Mathematics.}
Nauka, Moscow, 1983 (in Russian).

\bibitem{gj}
Knight's Tour Notes. 
Compiled by George Jelliss. 
Retrieved July 6, 2015.
http://www.mayhematics.com/t/t.htm

\bibitem{wp}
Knight's tour. In Wikipedia. Retrieved July 6, 2015, 
from the internet address https://en.wikipedia.org/wiki/Knight's\_tour.
\bibitem{biggs}
N. L. Biggs, E. K. Lloyd, R. J. Wilson. {Graph theory. 1736--1936.}
Second edition. The Clarendon Press, Oxford University Press, New York, 1986. 
\bibitem{wegener}
I. Wegener.
{\it Branching programs and binary decision diagrams. 
Theory and applications.} SIAM Monographs on Discrete Mathematics and Applications. 
SIAM, Philadelphia, PA, 2000. 
\bibitem{mk}
McKay, Brendan D.
Knight's Tours of an $8\times  8$ Chessboard.
Technical Report TR-CS-97-03 
(Department of Computer Science, Australian National University), 1997.

\bibitem{p}
Pickover, Clifford A.
The math book: From Pythagoras to the 57th dimension, 250 milestones in the history
of mathematics.
Sterling, New York, 
2009
\bibitem{seq}
The online encyclopedia of integer sequences. Founded by N.J.A. Sloane.
Retrieved July 6, 2015.
https://oeis.org/A165134
\bibitem{w}
Watkins, John J.  
Across the board: the mathematics of chessboard problems.
Princeton University Press, Princeton, NJ,  2004.
\end{thebibliography}
\end{document}